
\nopagenumbers
\footnote{{}}{\hskip -2 mm 2000 {\it Mathematics Subject
Classification}. Primary 11M06, Secondary 11F72.}

\def\txt#1{{\textstyle{#1}}}
\baselineskip=13pt
\def\hf{{\textstyle{1\over2}}}
\def\a{\alpha}
\def\d{{\,\rm d}}
\def\e{\varepsilon}
\def\f{\varphi}
\def\G{\Gamma}
\def\k{\kappa}
\def\s{\sigma}

\def\={\;=\;}
\def\zx{\zeta(\hf+ix)}
\def\zt{\zeta(\hf+it)}

\def\D{\Delta}
\def\no{\noindent}  \def\kk{\kappa_h}
\def\R{\Re{\rm e}\,} \def\I{\Im{\rm m}\,} \def\s{\sigma}
\def\z{\zeta}
\def\M{{\cal M}}

\def\no{\noindent} 

\def\H{H_j^3({\txt{1\over2}})}  \def\={\,=\,}
\def\hf{{\textstyle{1\over2}}}
\def\txt#1{{\textstyle{#1}}}
\def\f{\varphi}
\def\z#1{|\zeta (\hf +i#1)|^4}
\def\Z{{\cal Z}}
\font\tenmsb=msbm10
\font\sevenmsb=msbm7
\font\fivemsb=msbm5
\newfam\msbfam
\textfont\msbfam=\tenmsb
\scriptfont\msbfam=\sevenmsb
\scriptscriptfont\msbfam=\fivemsb
\def\Bbb#1{{\fam\msbfam #1}}

\def \NN {\Bbb N}
\def \CC {\Bbb C}
\def \RR {\Bbb R}
\def \ZZ {\Bbb Z}

\font\aa=cmss12 at 12pt
\font\bb=cmcsc10
\font\cc=cmcsc10 at 8pt

\font\ff=cmr9
\font\hh=cmbx12
\font\kk=cmcsc10 at 12pt
\def\rightheadline{\hfil\ff  Some    results for the
Riemann zeta-function and Hecke series\hfil\folio}
\def\leftheadline{\ff\folio\hfil Aleksandar Ivi\'c
 \hfil}
\def\emptyheadline{\hfil}
\headline{\ifnum\pageno=1 \emptyheadline\else
\ifodd\pageno \rightheadline \else \leftheadline\fi\fi}
\topglue2cm
\centerline{\aa  ON SOME  CONJECTURES AND RESULTS FOR THE  }
\medskip\centerline{\aa RIEMANN ZETA-FUNCTION
AND HECKE SERIES}
\bigskip\bigskip
\centerline{\kk Aleksandar Ivi\'c}
\bigskip\bigskip
\noindent {\bb Abstract.} {\ff We investigate the pointwise and mean square
order of the function $\Z_2(s)$, where
$\Z_k(s) = \int_1^\infty |\zx|^{2k}x^{-s}\d x, k\in\NN$. Three
conjectures involving $\Z_2(s)$ and certain exponential
sums of Hecke series in short intervals are formulated, which
have important consequences in zeta-function theory. A new order result
for $\Z_2(s)$ is obtained, and the function $\Z_k(s)$ is discussed.}

\bigskip\bigskip
{\hh 1. Spectral theory and the function $\Z_2(s)$}
\bigskip

\noindent

The spectral theory of the hyperbolic Laplacian has become increasingly
important in the theory of the Riemann zeta-function $\zeta(s)$,
especially in problems connected with the fourth moment of $|\zx|$.
For a comprehensive account of this subject we refer the reader
to Y. Motohashi's monograph [15]. In this section we shall briefly
state some necessary facts from spectral theory and introduce the
function $\Z_2(s)$, closely related to $|\zx|^4$.

\smallskip
Let $\,\{\lambda_j = \kappa_j^2 + {1\over4}\} \,\cup\, \{0\}\,$ be the
discrete spectrum of the hyperbolic Laplacian
$$
\Delta=-y^2\left({\left({\partial\over\partial x}\right)}^2 +
{\left({\partial\over\partial y}\right)}^2\right)
$$
acting over the Hilbert space composed of all
$\Gamma$-automorphic functions which are square integrable with
respect to the hyperbolic measure, where
$$
\G \;\cong\; SL(2,\,\ZZ)/\{+1,-1\}.
$$
 Let $\{\psi_j\}$ be a maximal
orthonormal system such that
$\Delta\psi_j=\lambda_j\psi_j$ for each $j\ge1$ and
$T(n)\psi_j=t_j(n)\psi_j$  for each integer $n\in\NN$, where
$$
\bigl(T(n)f\bigr)(z)\;=\;{1\over\sqrt{n}}\sum_{ad=n}
\,\sum_{b=1}^df\left({az+b\over d}\right)
$$
is the Hecke operator. We shall further assume that
$\psi_j(-\bar{z})=\e_j\psi_j(z)$ with $\e_j=\pm1$. We
then define ($s = \s + it$ will denote a complex variable)
$$
H_j(s)\;=\;\sum_{n=1}^\infty t_j(n)n^{-s}\qquad(\s > 1),
$$
which we call the Hecke series associated with the Maass wave form
$\psi_j(z)$, and which can be continued to an entire function.
As usual  we put
$$
\a_j = |\rho_j(1)|^2(\cosh\pi\kappa_j)^{-1},
$$
where $\rho_j(1)$ is the first Fourier coefficient of
$\psi_j(z)$. We note (see [15]) that
$$
\sum_{\k_j\le K}\a_jH_j^3(\hf) \ll K^2\log^CK \qquad(C > 0),\leqno(1.1)
$$
that $H_j(\hf) \ge 0$ (see Katok--Sarnak [12]), and that there are $
\asymp K$ eigenvalues $\k_j$ in $\,[K-1,\,K+1]\,$.

\medskip
The function $\Z_2(s)$ was introduced by Y. Motohashi
[14] (see also [10] and [15]), who showed  that  it
has meromorphic continuation over $\CC$. It is defined as
$$
\Z_2(s) \= \int_1^\infty |\zx|^4x^{-s}\d x \qquad(\s > 1).
$$
He also established that
in the half-plane $\s = \Re{\rm e}\, s > 0$
 it has the following singularities: the pole $s = 1$ of order five,
simple poles at $s = {1\over2} \pm i\k_j\,(\k_j
= \sqrt{\lambda_j - {1\over4}})$ and poles at $s = \hf\rho$,
where $\rho$ denotes complex
zeros of $\zeta(s)$. The residue of $\Z_2(s)$ at
$s = {1\over2} + i\k_h$ equals
$$
R(\k_h) = \sqrt{{\pi\over2}}{\Bigl(2^{-i\k_h}{\G({1\over4} - {i\over2}\k_h)
\over\G({1\over4} + {i\over2}\k_h)}\Bigr)}^3\G(2i\k_h)\cosh(\pi\k_h)
\sum_{\k_j=\k_h}\alpha_j H_j^3(\hf), \leqno(1.2)
$$\no
and the residue at $s = {1\over2} - i\k_h$ equals $\overline{R(\k_h)}$.
Thus there is an intrinsic connection between $\Z_2(s)$
and $\alpha_j H_j^3({1\over2})$, and therefore also between the fourth
moment of $|\zx|$ and $\alpha_j H_j^3({1\over2})$.
The plan of the paper is as follows: in Section 2 the
modified Mellin transform is studied and five lemmas needed in the
sequel are proved. Mean square results on $\Z_2(s)$ are studied in
Section 3, and the two conjectures on $\Z_2(s)$ and their
corollaries are given in Section 4. In Section 5 the conjecture on
exponential sums with $\a_j\H$ is made. A new pointwise bound for
$\Z_2(s)$ is obtained in Section 6. Finally in Section 7 the
general function $\Z_k(s)$ is studied.

\vfill\break
\topskip1cm
\bigskip
{\hh 2. The modified Mellin transform}
\bigskip

The function $\Z_2(s)$ represents a sort of a modified Mellin
transform of $|\zx|^4$, where the Mellin transform of an
integrable function $f(x)$ is commonly defined as
$$
\M[f(x)] \, = \, F(s) = \int_0^\infty x^{s-1}f(x)\d x\quad(s = \sigma + it).
$$
Mellin transforms play an important r\^ole in Analytic Number Theory.
By a change of variable they can be viewed as special cases of
Fourier transforms, and their theory built by using the theory of Fourier
transforms, for which the reader is referred to E.C. Titchmarsh's
classical monograph [16]. Namely
if $F(s) = {\cal M}[f(x)]$, then ($\xi = e^x$)
$$
F(\s+it) = \int_0^\infty \xi^{\s+it-1}f(\xi)\d \xi =
\int_{-\infty}^\infty e^{ixt}f(e^x)e^{\s x}\ dx
$$
is the Fourier transform of $f(e^x)e^{\s x}$.

An important   feature of Mellin transforms is
the so-called inversion formula. It states that
if $F(s) = \M[f(x)],\;$
$ y^{\sigma-1}f(y) \in L^1(0,\,\infty)$ and $f(y)$
is of bounded variation in a
neighbourhood of $y = x$, then
$$
{f(x + 0) + f(x - 0)\over2} \,
= \, {1\over2\pi i}\int\limits_{(\sigma)}\,F(s)x^{-s}\d s =
{1\over2\pi i}\lim_{T\to\infty}\int\limits_{\s-iT}^{\s+iT}\,F(s)x^{-s}\d s.
\leqno(2.1)
$$
We recall that
if $f(x)$ denotes measurable functions, then
$$
L^p(a,b) := \left\{f(x)\Big|\;
\int_a^b|f(x)|^p\d x < +\infty\right\}.
$$
The Mellin inversion formula is usually derived from the inversion
formula for Fourier transforms. Namely, if
$$
\widehat{g}(y) \, = \, \int_{-\infty}^\infty e^{ixy}g(x)\d x
$$
is the Fourier transform of $g(x)$, then under suitable conditions
this is equivalent to
$$
 g(x) \, = \, {1\over2\pi}\int_{-\infty}^\infty e^{-ixy}
\widehat{g}(y)\d y.    \leqno(2.2)
$$
For example
if $g,\,{\hat g} \in L^1(-\infty,\,\infty)$, then (2.2) holds
for almost all $x \in (-\infty,\,\infty)$. If additionally $g$
is continuous in $(-\infty,\,\infty)$, then (2.2) holds for
all $x\in(-\infty,\,\infty)$. A variant of
Parseval's formula for Fourier transforms is the identity
$$
{1\over2\pi}\int_{-\infty}^\infty|{\hat f}(x)|^2\d x =
\int_{-\infty}^\infty|f(x)|^2\d x, \leqno(2.3)
$$
and it can be used to derive Parseval's formula for Mellin transforms.

The modified Mellin transform $m[f(x)]$,
of which $\Z_2(s)$ is a special case, will be now defined as
$$
F^*(s) \= m[f(x)] \= \int_1^\infty f(x)x^{-s}\d x\qquad(s = \s+it),
\leqno(2.4)
$$
which is often more convenient to have in applications than the
ordinary Mellin transform, since now there are no convergence
problems at $x = 0$. If ${\bar f}(x) = f(1/x)$ when $0 < x \le 1$ and
${\bar f}(x) = 0$ otherwise, then
$$
m[f(x)] \= \M\left[{1\over x}{\bar f}(x)\right],\leqno(2.5)
$$
so that the properties of $m[f(x)]$ can be deduced from the properties
of the ordinary Mellin transform $\M[f(x)]$ by the use of (2.5).
In what  follows five lemmas concerning the properties
of the modified Mellin transform will be proved.

\medskip
LEMMA 1. {\it If $x^{-\s}f(x) \in L^1(1,\infty)$ and $f(x)$ is
continuous for $x > 1$, then}
$$
f(x) = {1\over2\pi i}\int_{(\s)}F^*(s)x^{s-1}\d s,\quad F^*(s) \=
m[f(x)]. \leqno(2.6)
$$

{\bf Proof.}  If $F^*(s) = m[f(x)]$, then from (2.5) and (2.1) we have
$$
{1\over x}{\bar f}(x) =  {1\over2\pi i}\int_{(\s)}F^*(s)x^{-s}\d s
\leqno(2.7)
$$
provided that $ y^{\sigma-1}{1\over y}{\bar f}(y) \in L^1(0,\,\infty)$.
This means that
$$
\int_0^1 y^{\sigma-1}{1\over y}\left|f\left({1\over y}\right)\right|
\d y  < +\infty,
$$
and making the change of variable  $y = 1/x$, we obtain that the
last condition is equivalent to $x^{-\s}f(x) \in L^1(1,\infty)$.
Changing $x$ to $1/x$ in (2.7) we obtain then (2.6).

\smallskip
LEMMA 2. {\it If $F^*(s) = m[f(x)],\,G^*(s) = m[g(x)]$, and $f(x), g(x)$
are real-valued, continuous functions for $x > 1$, such that}
$$
x^{{1\over2}-c}f(x) \in L^2(1,\infty),\quad x^{c-{1\over2}-\s}g(x)
\in L^2(1,\infty),
$$
{\it then}
$$
m[f(x)g(x)] = {1\over2\pi i}\int_{(c)}F^*(w)G^*(s + 1 - w)\d w.\leqno(2.8)
$$

{\bf Proof.} From   [16, Theorem 73] we have
$$
\int_0^\infty f(x)g(x)x^{s-1}\d x \= {1\over2\pi i}\int_{(c)}
F(w)G(s-w)\d w \leqno(2.9)
$$
if
$$\eqalign{&
F(s) = \M[f(x)],\; G(s) = \M[g(x)],\cr& x^{c-{1\over2}}f(x)
\in L^2(0,\,\infty),\;
x^{\s-c-{1\over2}}g(x) \in L^2(0,\,\infty).\cr}
$$
In place of $f(x)$ and $g(x)$ in (2.9) we shall take ${1\over x}\bar{f}(x)$
and $\bar{g}(x)$, respectively. By (2.5) we have $\M[{1\over x}\bar{f}(x)]
= F^*(s)$ and
$$
\M[\bar{g}(x)] \= \int_0^1 g\left({1\over x}\right)x^{s-1}\d x  \= G^*(s+1).
$$
Consequently (2.9) gives
$$
\int_0^\infty {1\over x}\bar{f}(x)\bar{g}(x)x^{s-1}\d x
= {1\over2\pi i}\int_{(c)} F^*(w)G^*(s+1-w)\d w.
$$
After the change of variable $y = 1/x$ this reduces to (2.8) if
$$
x^{c-{1\over2}}\cdot{1\over x}\bar{f}(x) \in L^2(0,\,\infty),
\quad
x^{\s-c-{1\over2}}\cdot{1\over x}\bar{g}(x) \in L^2(0,\,\infty),
$$
and these conditions are easily seen to be equivalent to
$$
x^{{1\over2}-c}f(x) \in L^2(1,\infty),\qquad
x^{c-{1\over2}-\s}g(x) \in L^2(1,\infty).
$$

\smallskip
LEMMA 3. {\it If $F^*(s) = m[f(x)],\,G^*(s) = m[g(x)]$, and $f(x), g(x)$ are
real-valued, continuous functions for $x > 1$, such that}
$$
x^{{1\over2}-\s}f(x) \in L^2(1,\infty),\quad x^{{1\over2}-\s}g(x)
\in L^2(1,\infty),
$$
{\it then}
$$
\int_1^\infty f(x)g(x)x^{1-2\s}\d x =  {1\over2\pi i}\int_{(\s)}
F^*(s)\overline{G^*(s)}\d s. \leqno(2.10)
$$

{\bf Proof.} Follows similarly to the  preceding proof from
Parseval's formula for Mellin transforms  in the form
$$
\int_0^\infty f(x)g(x)x^{2\s-1}\d x = {1\over2\pi i}\int_{(\s)}
F(s)\overline{G(s)}\d s
$$
if
$$\eqalign{
&F(s) = \M[f(x)],\; G(s) = \M[g(x)],\cr& x^{\s-{1\over2}}f(x)
\in L^2(0,\,\infty),\;x^{\s-{1\over2}}g(x)
\in L^2(0,\,\infty).\cr}
$$
The last relation follows e.g. from  [16, Theorem 72] by a
change of variable.

\medskip\no
Lemma 3, in the special case $f(x) = g(x)$,
is a natural tool for investigating mean square formulas connected with
the modified Mellin transform. In particular, it offers the possibility
to obtain mean square estimates of $f(x)$ from the mean square estimates
of $m[f(x)]$, provided we have adequate analytic information about the
latter. A result in this direction, which is useful for the applications
that we have mind, will be given now as

\medskip
LEMMA 4. {\it Suppose that $g(x)$ is a real-valued,
integrable function on $[a,b]$, a subinterval
of $[2,\,\infty)$, which is not necessarily finite. Then}
$$
\int\limits_0^{T}\left|\int\limits_a^b g(x)x^{-s}\d x\right|^2\d t
\le 2\pi\int\limits_a^b g^2(x)x^{1-2\s}\d x \quad(s = \s+it\,,T > 0,\,a<b).
\leqno(2.11)
$$
\medskip
{\bf Proof.} Let
$$
I := \int_0^{T}\left|\int_a^b g(x)x^{-s}\d x\right|^2\d t
\qquad(s = \s + it,\;T  > 0).\leqno(2.12)
$$
In Lemma 3 set $f(x) = g(x)$
if $a \le x \le b$ and $f(x) = 0 $ otherwise. Then
$$
F^*(s) \= m[f(x)] \= \int_a^b g(x)x^{-s}\d x,\quad F^*(s) \= G^*(s),
$$
and $x^{{1\over2}-\s}f(x) \in L^2(1,\,\infty)$ for any $\s$.
Consequently (2.10) of Lemma 3 (with $f \equiv g$) gives
$$
{I\over2\pi} \;\le\; {1\over2\pi}\int_{-\infty}^\infty |F^*(\s + it)|^2\d t
= \int_a^b g^2(x)x^{1-2\s}\d x,\leqno(2.13)
$$
and (2.11) follows from (2.12) and (2.13).

\medskip
In case when we have sufficient information on $f$ and $F^*$, we
can obtain an explicit bound for the mean square estimate of $f(x)$.
Such a result is furnished by

\medskip
LEMMA 5. {\it Let $f(x) \in C^\infty[2,\infty)$ be a real-valued function
for which}

\vskip0.5mm

a) $\int_1^X|f^{(r)}(x)|\d x \ll_{\e,r} X^{1+\e}\;(r = 0,1,2,\,\ldots\,)$
{\it and}

\vskip0.5mm

b) $F^*(s) = m[f(x)]$ {\it has a pole at $s =1 $ of order $\ell$,
 but otherwise can be
analytically continued to the region $\R s > \hf$, where it is of
polynomial growth in $|\I s|$. Then for
$\hf < \s < 1$ and any given $\e >0$ we have}
$$
\int_T^{2T}f^2(x)\d x \ll_\e \log^{\ell-1}T\cdot\int_{T/2}^{5T/2}
|f(x)|\d x +
T^{2\s-1}\int_0^{T^{1+\e}}|F^*(\s + it)|^2\d t.
\leqno(2.14)
$$

\medskip
{\bf Proof.} Let $\f(x) \in C^\infty(0,\infty)$
be a test function such  that $\f (x) \ge0,\,\f(x) = 1$ for
$T \le x\le 2T$, $\f (x) = 0$ for $x < \hf T$ or $x > {\txt{5\over2}}T$
$(T \ge T_0> 0)$, $\f(x)$ is increasing in $[\hf T,T]$ and
decreasing in $[2T,{\txt{5\over2}}T]$. Then obviously
$$
I_1 \;\le\;I_2,
$$
where
$$
I_1 \;:=\; \int_T^{2T}f^2(x)\d x,\qquad
I_2 \;:=\; \int_{T/2}^{5T/2}\f(x)f^2(x)\d x.
$$
From the assumption a) we have $x^{-\s}f(x)\in L^1(1,\,\infty)$ if
$\s > 1$. Therefore Lemma 1 gives
$$
f(x) = {1\over2\pi i}\int_{(1+\e)}F^*(w)x^{w-1}\d w
= {1\over2\pi i}\int_{{\cal L}}F^*(w)x^{w-1}\d w + Q_{\ell-1}(\log x),
\leqno(2.15)
$$
where $Q_{\ell-1}(\log x)$ is a polynomial in $x$ of degree $\ell-1$,
and ${\cal L}$ is the line $\R w = 1+\e$ with a small indentation to
the left at the pole $w=1$ of $F^*(w)$ of order $\ell$, so that
by the residue theorem we pick a contribution equal to
$Q_{\ell-1}(\log x)$. Therefore (2.15) yields
$$\eqalign{
I_2 &=
{1\over2\pi i}\int_{{\cal L}}F^*(w)\left(
\int_{T/2}^{5T/2}\f(x)f(x)x^{w-1}\d x\right)\d w \cr&+
O\left(\log^{\ell-1}T\int_{T/2}^{5T/2}\f(x)|f(x)|\d x \right).\cr}
\leqno(2.16)
$$
We integrate $r$ times by parts to
obtain
$$\eqalign{
&\int_{T/2}^{5T/2}\f(x)f(x)x^{w-1}\d x\cr&
= {(-1)^r\over w\cdots(w+r-1)}
\int_{T/2}^{5T/2}(\f(x)f(x))^{(r)}x^{w+r-1}\d x.\cr}\leqno(2.17)
$$
Since $\f^{(r)}(x) \ll_r T^{-r}$ and a) of the lemma holds, then
by (2.17) it follows that the portion of the integral over $w$ in (2.16)
for which $|v| > T^{1+\e}$ makes
a negligible contribution, namely $\ll T^{-A}$ for any given $A > 0$,
provided that we choose $r = r(A, \e)$ sufficiently large. This
remains true
even if we move the contour of integration to the left. Hence if
$0 < u < \hf,\, w = u + iv$, then by (2.16) and the Cauchy-Schwarz
inequality for integrals we infer that
$$
\eqalign{
&I_2 =
 \,{1\over2\pi}\int_{-T^{1+\e}}^{T^{1+\e}}
F^*(u+iv)\int_{T/2}^{5T/2}\f(x)f(x)x^{u+iv-1}\d x\d v \cr&
\quad+  O\left(\log^{\ell-1}T\cdot\int_{T/2}^{5T/2}\f(x)|f(x)|\d x \right)
\cr&  \,\ll\;\left\{\int_{0}^{T^{1+\e}}|F^*(u+iv)|^2\d v
\int\limits_{0}^{T^{1+\e}}\left|
\int\limits_{T/2}^{5T/2}\f(x)f(x)x^{u+iv-1}\d x\right|^2\d v\right\}^{1/2}
\cr&
\;+ \;\log^{\ell-1}T\cdot\int_{T/2}^{5T/2}\f(x)|f(x)|\d x.
\cr}
$$
Now we apply Lemma 4 (with $g(x) = \f(x)f(x),\,s= 1 - u +iv$, $T$
replaced by $T^{1+\e}$) to deduce that
$$
\eqalign{&
\int_{0}^{T^{1+\e}}\left|
\int_{T/2}^{5T/2}\f(x)f(x)x^{u+iv-1}\d x\right|^2\d v
\cr&
\ll \int_{T/2}^{5T/2}\f^2(x)f^2(x)x^{2u-1}\d x
\ll T^{2u-1}\int_{T/2}^{5T/2}\f(x)f^2(x)\d x = T^{2u-1}I_2,
\cr}
$$
since $0 \le \f(x) \le 1$. It follows that
$$
I_2 \ll \log^{\ell-1}T\cdot\int\limits_{T/2}^{5T/2}\f(x)|f(x)| \d x
+ \left(\int\limits_{0}^{T^{1+\e}}
|F^*(u+iv)|^2\d v\cdot T^{2u-1}I_2\right)^{1/2},
$$
which easily gives (2.11) with $\s = u$.

\bigskip
{\hh 3. The mean square of $\Z_2(s)$}
\bigskip

Mean square problems involving $\Z_2(s)$ are naturally of interest.
They were investigated in [10], where it was shown that
$$
\int_0^T|\Z_2(\s + it)|^2\d t \;\ll\;T^{10-8\s\over3}\log^CT
\qquad(\hf < \s <1,\,C > 0).\leqno(3.1)
$$
We introduce now the function $E_2(T)$, the error term in the asymptotic
formula for the mean fourth power of $\zt$, customarily
defined by the relation
$$
\int_0^T|\zt|^4\d t \;=\; TP_4(\log T) \;+\;E_2(T),\leqno(3.2)
$$
with
$$
P_4(x) \;=\; \sum_{j=0}^4\,a_jx^j, \quad a_4 = {1\over2\pi^2}.\leqno(3.3)
$$
For some recent results on $E_2(T)$ see [3]-[6], [8], [9], [14] and [15].
For the explicit evaluation of the $a_j$'s in (3.3), see [4].
A fundamental result in the theory of $E_2(T)$ is the mean square estimate
$$
\int_0^T E_2^2(t)\d t \;\ll_\e\; T^{2+\e},\leqno(3.4)
$$
where $\e$ denotes arbitrarily small positive constants, not
necessarily the same ones at each occurrence. In fact, Ivi\'c--Motohashi
[8] proved (3.4) with $\log^CT$ in place of $T^\e$, but for our
present purpose (3.4) suffices; the integral in (3.4) is actually
$\gg T^2$, as shown in [6], so that the bound in (3.4) is
essentially best possible. We have (see  [5], [6], [9] and [15])
$$
E_2(T) \;\ll\;T^{2/3}\log^8T,\quad E_2(T) \= \Omega_\pm(\sqrt{T})
\leqno(3.5)
$$
and
$$
\int_0^T E_2(t)\d t \= O(T^{3/2}),\quad  \int_0^T E_2(t)\d t
\= \Omega_\pm(T^{3/2}).
\leqno(3.6)
$$
Thus if $c$ is such a constant for which
$$
E_2(T) \;\ll_\e\;T^{c+\e}\leqno(3.7)
$$
holds, then (3.5) implies that one must have
$\hf \le c \le {2\over3}$. It was proved in [10]
that, besides (3.1), one also has
$$
\int_0^T|\Z_2(\s + it)|^2\d t \;\ll_\e\;
T^\e\left(T + T^{2-2\s\over1-c}\right)
\qquad(\hf < \s < 1).\leqno(3.8)
$$

We can derive a mean square result on $\Z_2(s)$ by the use of Lemma 3.
Namely from (3.2) we obtain, by integration by parts, the
representation
$$\eqalign{
\Z_2(s) &\= \int_1^\infty \left(P_4(\log x) + P_4'(\log x)
+ E_2'(x)\right)x^{-s}\d x\cr&
\= \sum_{j=0}^5c_j(s-1)^{-j} +
s\int_1^\infty E_2(x)x^{-s-1}\d x.\cr}\leqno(3.9)
$$
It follows from (3.4) and the Cauchy-Schwarz inequality for
integrals that the function
$$
z_2(s) := {\cal Z}_2(s) - \sum_{j=0}^5c_j(s-1)^{-j} =
s\int_1^\infty E_2(x)x^{-s-1}\d x\leqno(3.10)
$$
is regular for $\s > \hf$, where the constants $c_j$ may be
explicitly evaluated in term of the $a_j$'s. Then Lemma 3 gives
$$
\int_1^\infty E_2^2(x)x^{-1-2\s}\d x =
 {1\over 2\pi}\int_{-\infty}^\infty \left|z_2(\s + it)
\over\s+it\right|^2\d t\quad(\s > \hf).\leqno(3.11)
$$
Thus, for $\hf < \s < 1$, we have
$$
1\gg \int_T^{2T}\left|z_2(\s + it)\over\s+it\right|^2\d t
\gg T^{-2}\int_T^{2T}\left|z_2(\s + it)\right|^2\d t,
$$
which yields
$$
\int_T^{2T}|{\cal Z}_2(\s + it)|^2\d t \;\ll\;
\int_T^{2T}\left(|z_2(\s + it)|^2 + 1\right)\d t  \;\ll\;
T^2\;(\hf < \s < 1).
\leqno(3.12)
$$
Note that if (3.12)  is known, then (3.11) yields (3.4).
To see this let $s = \sigma + it$ with
$\hf < \sigma < 1, \,t \ge t_0 > 0$. Then by the residue theorem
$$
{\cal Z}_2(s) = {1\over2\pi i}\int_{\cal D}X^w\G(w){\cal Z}_2(w+s)\d w
\qquad(2 \le X \ll t^A),\leqno(3.13)
$$
where ${\cal D}$ is the rectangle with vertices $\hf - \sigma + \e
\pm i\log^2t, \,1-\sigma +\e \pm i\log^2t$ and $0 < \e < \sigma - \hf$.
By Stirling's formula  for the gamma-function  it follows then that
$$
{\cal Z}_2(s) \ll_\e X^{\hf-\sigma+\e}\int_{-\log^2t}^{\log^2t}e^{-|v|}
|{\cal Z}_2(\hf+\e+iv+it)|\d v + X^{1-\sigma+\e},\leqno(3.14)
$$
since ${\cal Z}_2(s) \ll_\s 1$ for $\s > 1$. If we additionally
suppose that $T \le t\le 2T$, then from (3.12) (with $\sigma = \hf+\e$)
and (3.14) we obtain
$$\eqalign{&
\int_T^{2T}|{\cal Z}_2(\sigma + it)|^2\d t\cr&
\ll_\e TX^{2-2\sigma+\e} + X^{1-2\sigma+\e}\int_{-\log^2T}^{\log^2T}e^{-|v|}
\int_T^{2T}|{\cal Z}_2(\hf + \e  +iv + it)|^2\d t\,\d v\cr&
\ll_\e TX^{2-2\sigma+\e} + T^2X^{1-2\sigma+\e} \ll_\e T^{3-2\sigma+\e}\cr}
$$
if we choose $X = T$. Hence it follows that
$$
\int_1^T|{\cal Z}_2(\sigma + it)|^2\d t \;\ll_\e\;
T^{3-2\sigma+\e}\qquad(\hf < \sigma < 1).\leqno(3.15)
$$
Since $3 - 2\sigma < 2$ for $\sigma > \hf$, (3.15) yields
$$
\int_{-\infty}^\infty \left|{z_2(\sigma+it)\over\sigma+it}
\right|^2\d t \ll 1 \qquad(\hf < \sigma < 1),
$$
and consequently (3.11) yields (3.4) with $\sigma = {1+\e\over2}$.
As proved in [8],  (3.4)  easily gives then
$$
E_2(T) \ll_\e T^{{2\over3}+\e},\;\int_0^T |\zt|^{12}\d t \ll_\e T^{2+\e},
\leqno(3.16)
$$
and both of these estimates are best known up to ``$\e$".
This shows the strength and importance of mean square estimates for
${\cal Z}_2(s)$ (see also (3.29)).

\smallskip

A natural problem is to ask: what is the true order of magnitude of
$$
I_\s(T) \;:=\; \int_1^T|\Z_2(\s+it)|^2\d t \qquad(\s > \hf)?\leqno(3.17)
$$
It seems very hard to speculate what ought to be the shape of the
asymptotic formula for $I_\s(T)$ for fixed $\s > \hf$, not to
mention values like $\s = \hf + {1\over\log T}$
(in view of the poles $s = \hf\pm i\k_j$ of $\Z_2(s)$)
etc. The lower limit of integration in (3.17) is unity and not
zero for technical reasons, in view of the pole $s=1$ of $\Z_2(s)$.
The bound
$$
I_{{1\over2}+\e}(T) \;\ll\;T^2 \leqno(3.18)
$$
follows from (3.15)
and in view of the fact that (3.4) is essentially best possible,
it follows that by using Parseval's formula, namely
the identity (3.11), we
cannot obtain a stronger estimate than (3.18). Nevertheless
there seems to be no apparent reason why (3.18) could not be improved.
It can be also remarked that (3.18), via the second bound in
(3.16) and the Cauchy-Schwarz
inequality for integrals, yields the bound
$$
\int_0^T|\zt|^8\d t \;\ll_\e\; T^{{3\over2}+\e},\leqno(3.19)
$$
which is  best known up to ``$\e$". Thus any improvement of
(3.18) would have important consequences in zeta-function theory.
On the other hand it is of interest to obtain lower bounds (i.e.,
omega-results) for $I_\s(T)$. In this direction we have

\medskip
{\bf THEOREM 1}. {\it For any given $\e > 0$ and fixed $\s$ satisfying
$\hf < \s < 1$ we have}
$$
\int_1^T|\Z_2(\s + it)|^2\d t \;\gg_\e\; T^{2-2\s-\e}.\leqno(3.20)
$$

\medskip\no {\bf Proof.} Let, as usual,
$$
Z(t) := \chi^{-{1\over2}}(\hf+it)\zt,\quad \chi(s) = {\zeta(s)\over
\zeta(1-s)} = {\pi^{s-{1\over2}}\G(\hf-\hf s)\over\G(\hf s)}.
$$
Then $Z(t) \in \RR$ if $t\in \RR$, $|Z(t)| = |\zt|$, and $Z^4(t)$ and
all its derivatives are bounded in mean by a suitable log-power. This
follows by H\"older's inequality for integrals, and the fact that
(see [11, Chapter 3])
$$
\eqalign{& Z^{(k)}(t) \;=\; O_k\left(t^{-1/4}({\txt{3\over2}}\log t
)^{k+1}\right) \;+
\cr&
+ \;2\sum_{n\le\sqrt{{t\over2\pi}}}n^{-{1\over2}}
\Bigl(\log{\sqrt{t/(2\pi)}\over n}\Bigr)^k
\cos\left(t\log{\sqrt{t/(2\pi)}\over n}
- {t\over2} - {\pi\over8} + {\pi k\over2}\right).\cr}
$$
Another way to see this is to use Leibniz's formula for the derivative
of a product, the expression for $\chi(s)$ and properties of the
gamma-function.
If we set $f(x) = Z^4(x)$, then $F^*(s) = \Z_2(s)$, and the
assumptions a) and b) of Lemma 5 are satisfied ($\ell = 5,\,\Z_2(s)$
is of polynomial growth for $\s > \hf$ by (3.9)). Therefore (2.14) gives,
for $\hf < \s < 1$,
$$\eqalign{
\int_T^{2T}|\zt|^8\d t &\ll_\e \log^4T\cdot \int_{T/2}^{5T/2}|\zt|^4\d t\cr&
+ T^{2\s-1}\int_0^{T^{1+\e}}|\Z_2(\s+it)|^2\d t.\cr}\leqno(3.21)
$$
But (see [3, Theorem 6.5]) we have
$$
\int_T^{2T}|\zt|^8\d t \;\gg\; T\log^{16}T,
$$
hence (3.21) yields
$$
T\log^{16}T \;\ll_\e\; T^{2\s-1}\int_0^{T^{1+\e}}|\Z_2(\s+it)|^2\d t.
\leqno(3.22)
$$
The lower bound in (3.20) follows then if we substitute $T$ for $T^{1+\e}$
in (3.22).

\smallskip

The mean square of $\Z_2(s)$ possesses a convexity property, embodied in

\medskip
{\bf THEOREM 2}. {\it For fixed $c,\s$ such that $\hf < c \le \s$
and $t \ge T \ge 3$ we have}
$$\eqalign{
\int\limits_T^{2T}|\Z_2(\s+it)|^2\d t &\ll_\e T^{-1}\cr&
+ T^\e\left(T^{2c-2\s}\int\limits_T^{2T}|\Z_2(c+it)|^2\d t
+ T^{2-2\s} + T^{5-6\s}\right),\cr}\leqno(3.23)
$$
$$
\Z_2(\s+it) \ll_\e t^{-1} + t^{c-\s+\e}\int_{t-t^\e}^{t+t^\e}
|\Z_2(c+iv)|\d v + t^{3-4\s}\log^4t.
\leqno(3.24)
$$

\medskip\no
{\bf Proof.} The proof of both bounds is similar (for (3.24) we take
$X=t^4$), so only (3.23) will be considered in detail.
Clearly we may assume that $c <\s$. Let $T \le t \le 2T
\le X \le T^A$, where $A \,( > 1)$ is a constant. For $\s > 1$ we have
$$
\eqalign{
\Z_2(s) &= \int_1^{2T^{1-\e}}\rho(x)|\zx|^4x^{-s}\d x
+ \int_{T^{1-\e}}^{2X}\s(x)|\zx|^4x^{-s}\d x \cr&
+ \int_X^\infty\tau(x)|\zx|^4x^{-s}\d x = I_1(s)
+ I_2(s) + I_3(s),\cr}
$$
say. Here $\rho(x), \s(x), \tau(x) \in C^\infty$ are nonnegative
functions such that: $\rho(x) = 1$ for $1 \le x \le T^{1-\e}$,
$\rho(x) = 0$ for $x \ge 2T^{1-\e},\, \s(x) = 0$ for $x \ge 2X,
\tau(x) = 1 - \s(x)$. By  repeated integration by parts,
similarly as in the proof of (3.21), we obtain
$$
I_1(s) \= {|\zeta(\hf+i)|^4\over s-1} + O(t^{-2}) \ll {1\over T}.
$$
Next note that
$$
\eqalign{
I_3(s) &= \int_X^\infty \tau(x)Q_4(\log x)x^{-s}\d x
+ \int_X^\infty \tau(x)E_2'(x)x^{-s}\d x\cr&
= {1\over s-1}\int_X^\infty x^{1-s}\Bigl(\tau'(x)Q_4(\log x)
+ \tau(x)x^{-1}Q_4'(\log x)\Bigr)\d x\cr&
- \int_X^\infty E_2(x)\Bigl(\tau'(x)x^{-s}-s\tau(x)x^{-s-1}\Bigr)\d x
= I_4(s) - I_5(s),\cr}
$$
say. Since $\tau'(x) \ll X^{-1}$ and $\tau'(x) = 0$ for $x \ge 2X$,
it follows that
$$
\eqalign{
I_4(s) &
= {1\over s-1}\int_X^{2X}x^{1-s}\Bigl(\tau'(x)Q_4(\log x)
+ \tau(x)x^{-1}Q_4'(\log x)\Bigr)\d x\cr&
+ {1\over s-1}\int_{2X}^\infty x^{-s}Q_4'(\log x)\d x
= I_6(s) + I_7(s),\cr}
$$
say. The function $I_6(s)$ is regular for $s \not= 1$ and for
$T \le t \le 2T$ it is $\ll T^{-1}X^{1-\s}\log^4X$. We can
evaluate $I_7(s)$ and obtain
$$
I_7(s) \= {c(2X)^{1-s}\over(s - 1)^4},
$$
which provides analytic continuation of $I_7(s)$ to $\CC$.
In view of (3.4) (or (3.6)) $I_5(s)$ is regular for $\s > \hf$,
and therefore we obtain
$$
\int_T^{2T}|I_5(s)|^2\d t \ll T^2\int_X^\infty E_2^2(x)x^{-1-2\s}\d x
\ll T^2X^{1-2\s}\log^CX    \leqno(3.25)
$$
by using Lemma 4. Consequently we have
from (3.24) and (3.25), for $\s > \hf$,
$$
\eqalign{
&\int_T^{2T}(|I_1(s)|^2+|I_3(s)|^2)\d t \cr&\ll T^{-1}
+ T^{-1}X^{2-2\s}\log^8X + T^2X^{1-2\s}\log^CX\cr&
\ll T^{-1} + T^{5-6\s}\log^CT\cr}
$$
for $X = T^3$, which we henceforth assume. It remains to deal with
the mean square of $I_2(s)$, which we write as a sum of $O(\log T)$
integrals of the type
$$
F_K(s) \;:=\; \int_{K/2}^{5K'/2}\f(x)|\zx|^4x^{-s}\d x
\quad(T^{1-\e} \le K < K' \le 2K \ll X),
$$
where $\f(x) \in C^\infty$ is  a nonnegative function supported
in $[K/2,\,5K'/2]$ such that $\f(x) = 1$ for $K < K' \le 2K$, and
$$
\f^{(r)}(x) \;\ll_r\;K^{-r}\qquad(r = 0,1,2,\ldots).\leqno(3.26)
$$
To connect $F_K(s)$ and $\Z_2(s)$ note that from
the Mellin  inversion formula (2.6) we have
$$
|\zx|^4 \= {1\over2\pi i}\int_{(1+\e)}{\cal Z}_2(s)x^{s-1}\d s \qquad(x>1).
$$
Here we replace the line of integration  by the contour ${\cal L}$,
consisting of the same straight line from which the segment
$[1+\e-i,\,1+\e+i]$ is removed and replaced by a circular arc
of unit radius, lying
to the left of the line, which passes over the pole $s =1 $ of
the integrand. By the residue theorem we have
$$
|\zx|^4 \= {1\over2\pi i}\int_{\cal L}{\cal Z}_2(s)x^{s-1}\d s
+ Q_4(\log x) \qquad(x > 1),\leqno(3.27)
$$
where we have set (cf. (3.2))
$Q_4(\log x) \= P_4(\log x) + P'_4(\log x)$. Hence by using (3.27)
we obtain
$$\eqalign{
F_K(s) &= {1\over2\pi i}\int_{\cal L}{\cal Z}_2(w)
\left(\int_{K/2}^{5K'/2}\f(x)x^{w-s-1}\d x\right)\d w\cr&
+ \int_{K/2}^{5K'/2}\f(x)Q_4(\log x)x^{-s}\d x.\cr}\leqno(3.28)
$$
In view of (3.26) we infer, by repeated integration by parts,
that the last integral in (3.28) is $\ll T^{-A}$ for any given $A>0$.
Similarly we note that
$$\eqalign{
&\int_{K/2}^{5K'/2}\f(x)x^{w-s-1}\d x\cr&
= (-1)^r\int_{K/2}^{5K'/2}\f^{(r)}(x){x^{w-s+r-1}\over(w-s)\cdots
(w-s+r-1)}\d x \ll T^{-A}\cr}
$$
for any given $A>0$, provided that $|\I w - \I s| > T^\e$ and
$r = r(A,\e)$ is sufficiently large. Thus if in the $w$--integral
in (3.28) we replace the contour ${\cal L}$ by the straight line
$\R w = c$, we shall obtain
$$
F_K(s) \;\ll\; K^{c-\s}\int_{t-T^\e}^{t+T^\e}|\Z_2(c+iv)|\d v + T^{-2},
$$
which gives
$$\eqalign{&
\int_T^{2T}|I_2(s)|^2\d t\cr&
\ll T^{-1} + \log T\max_{T^{1-\e}\ll K\ll T^3}\,K^{2c-2\s}\int_T^{2T}\left(
\int_{t-T^\e}^{t+T^\e}|\Z_2(c+iv)|\d v\right)^2\d t\cr&
\ll T^{-1} + \log T\max_{T^{1-\e}\ll K\ll T^3}\,K^{2c-2\s}T^{\e}\int_T^{2T}
\int_{t-T^\e}^{t+T^\e}|\Z_2(c+iv)|^2\d v\,\d t\cr&
\ll T^{-1} + \log T\max_{T^{1-\e}\ll K\ll T^3}\,K^{2c-2\s}T^\e
\int\limits_{T-T^\e}^{2T+T^\e}|\Z_2(c+iv)|^2
\int\limits_{v-T^\e}^{v+T^\e}\d t\,\d v\cr&
\ll T^{-1} + T^{2c-2\s+4\e}\int_{T}^{2T}|\Z_2(c+iv)|^2\d v + T^{2-2\s+4\e},
\cr}
$$
where we used the bound $\Z_2(\s+it) \ll_\e t^{1-\s+\e}$ (see Theorem 3),
whose proof is independent of the present theorem. Collecting the above
bounds we obtain (3.23).

\medskip
{\bf Corollary 1.} {\it For any given $\e > 0$ we have}
$$
\int_T^{2T}|\Z_2(\s+it)|^2\d t \ll_\e T^{-1} + T^{{8\over3}-2\s+\e}
\qquad(\s \ge {\txt{5\over6}}).\leqno(3.29)
$$

\medskip\no
{\bf Proof of Corollary 1.} From (3.8) (with $c = {2\over3}$) we have
$$
\int_T^{2T}|\Z_2(\s+it)|^2\d t \ll_\e T^{1+\e}\qquad({\txt{5\over6}}
\le \s \le 1),
$$
hence (3.29) follows from (3.23) (with $c = {5\over6}$) since
$5 -  6\s \le {8\over3}- 2\s$ for $\s \ge {7\over 12}$. Note that
(3.29) sharpens (3.8) (with the best known value $c = {2\over3}$)
for  $\s \ge {\txt{5\over6}}$.

\bigskip
{\hh 4. The conjectures on $\Z_2(s)$}
\bigskip\no

It seems reasonable that the lower bound (3.20) of Theorem 1
is essentially of the correct order of magnitude. Therefore we formulate
the following

\medskip
{\bf Conjecture 1.} {\it For any given $ \e > 0,
\hf < \s < 1$ and $T \gg 1$ we have}
$$
\int_1^T|\Z_2(\s+it)|^2\d t  \;\ll_{\e}\; T^{2-2\s+\e}. \leqno(4.1)
$$
This conjecture is very strong. It implies
the essentially best possible bounds
for $E_2(T)$ and the eighth moment of $|\zt|$. This is contained in

\smallskip {\bf Corollary 2.} {\it If Conjecture 1 holds, then}
$$
E_2(T) \;\ll_{\e}\; T^{{1\over2}+\e}.\leqno(4.2)
$$

\smallskip {\bf Corollary 3.} {\it If Conjecture 1 holds, then}
$$
\int_0^T|\zt|^8\d t \;\ll_{\e}\; T^{1+\e}.\leqno(4.3)
$$
\smallskip\no
{\bf Proof of Corollary 2}.
From the defining relation (3.2) it is not difficult to obtain
(see e.g., [5, (5.3)]) that $(C_1, C_2 > 0,\;1 \ll H \le {\txt{1\over4}}T)$,
$$
E_2(T) \le C_1H^{-1}\int_T^{T+H}E_2(x)f(x)\d x + C_2H\log^4T,\leqno(4.4)
$$
where $f(x)\;( > 0)\,$ is a smooth function supported in $\,[T,\,T+H]\,$,
such that $f(x) = 1$ for $T + {1\over4}H \le x \le T + {3\over4}H$.
If we integrate (3.27) from $x= 1$ to $x = T$ and
take into account the defining relation (3.2) of $E_2(T)$, we shall obtain
$$
E_2(T) \= {1\over2\pi i}\int_{\cal L}{\cal Z}_2(s){T^s\over s}\d s
+ O(1)\qquad(T > 1).\leqno(4.5)
$$
Then from (4.4) and (4.5) we have ($\hf < c < 1,\,T>1$)
$$
E_2(T) \le {C_1\over2\pi iH}\int_{(c)}{{\cal Z}_2(s)\over s}
\int_T^{T+H} f(x)x^s\d x\,\d s + C_2H\log^4T,\leqno(4.6)
$$
and we also have an analogous lower bound for $E_2(T)$.
Since $f^{(r)}(x) \ll_r H^{-r}$ it follows that the $s$--integral
in (4.6) can be truncated at $|\Im{\rm m}\, s|= T^{1+\e}H^{-1}$ with
a negligible error, for any $c$ satisfying $\hf < c < 1$.
We take $c = \hf + \e$ and use (4.1), coupled with the
Cauchy-Schwarz inequality for integrals, to deduce that
$$
\int_{({1\over2}+\e)}{{\cal Z}_2(s)\over s}
\int_T^{T+H} f(x)x^s\d x\,\d s \ll_\e H^{1/2}T^{1+2\e},
$$
so that Corollary 2 follows from (4.6) with $H = T^{1\over2}$.
\smallskip\no
{\bf Proof of Corollary 3}. From (3.21) we have
$$
\int_T^{2T}|\zt|^8\d t \ll_\e
T^{2\s-1}\int_0^{T^{1+\e}}|\Z_2(\s+it)|^2\d t \quad(\hf < \s < 1),
\leqno(4.7)
$$
so that (4.3) follows from (4.1) and (4.7).

It is plausible that Conjecture 1 is
equivalent to Corollary 3. We can prove something a little weaker.
Namely if (4.3) holds, then by the above method one can
sharpen (3.28) to
$$
\int_T^{2T}|\Z_2(\s + it)|^2\d t \;\ll_\e\;
\cases{T^{4-4\s+\e}\qquad(\hf < \s \le 1),\cr
\cr
T^{2-2\s+\e} + T^{-1}\quad(\s \ge 1).
\cr}\leqno(4.8)
$$
Conversely, (4.8) (with $\s = 1-\e$) implies (4.3) by (4.7), so that
(4.3) and (4.8) are equivalent.

\medskip
In view of the discussion on the true order of $I_\s(T)$,
it seems in place to discuss also the problem of the order of
$\Z_2(\s + it),\, t \ge t_0 > 0$ and $\hf < \s < 1$. Conjecture 1 says
that $\Z_2(\s + it)$ is small in mean square. Perhaps it is even also small
pointwise, so  the following conjecture is now proposed.

\medskip
{\bf Conjecture 2.} {\it For any given $\e > 0$ we have}
$$
\Z_2(\s+it) \;\ll_{\e}\; |t|^\e \qquad(\s > \hf). \leqno(4.9)
$$

\medskip

Conjecture 2 is the analogue of the
classical Lindel\"of hypothesis  ($\zt  \;\ll_{\e}\; |t|^\e$) in
the equivalent form
$$
\zeta(\s + it)  \;\ll_{\e}\; |t|^\e\qquad(\s > \hf).\leqno(4.10)
$$
Since both $\Z_2(s)$ and $\zeta(s)$ take conjugate
values at conjugate points, it suffices in (4.9) and (4.10) to assume
that $t > 0$. Conjecture 2 implies Conjecture 1, which easily follows from

\medskip
{\bf Corollary 4.} {\it Conjecture 2 is equivalent to the statement that,
for any given $\e > 0$ and $\hf < \s < 1$},
$$
\Z_2(\s + it) \;\ll_\e\; t^{{1\over2}-\s+\e}\qquad(t \ge t_0 > 0).
\leqno(4.11)
$$

\medskip\no
{\bf Proof of Corollary 4.}  Trivially (4.11) implies (4.9), so we have only
to prove that (4.9) implies (4.11).
We suppose that $\s > 1$ and proceed similarly as in the proof of (3.23)
to obtain
$$
\eqalign{
\Z_2(s) &= O({1\over t}) + {1\over2\pi i}\int_{{\cal L}}\Z_2(w)
\left(\int_{T^{1-\e}}^\infty(1-\rho(x))x^{w-s-1}\d x\right)\d w\cr&
+ \int_{T^{1-\e}}^\infty(1-\rho(x))Q_4(\log x)x^{-s}\d x
\quad(T\le t\le 2T).\cr}\leqno(4.12)
$$
An integration by parts shows that the last integral above is $\ll
t^{\e-1}$. Since $(1-\rho(x))^{(r)} \ll_r T^{r(\e-1)}$ for
$r = 0,1,2,\ldots\,$, it follows by $r$ integrations by parts that
$$
\int_{T^{1-\e}}^\infty(1-\rho(x))x^{w-s-1}\d x \;\ll\;T^{-A}
$$
for any given $A > 0$ if $\,|\I w - \I s| > T^\e\,$ and $r = r(\e,A)$
is sufficiently large. Hence (4.12) yields, on replacing ${\cal L}$ by
the line $\R w = \hf + \delta$,
$$\eqalign{
\Z_2(s) &\ll t^{\e-1} + t^\e\max_{|v-t|\le t^\e}|\Z_2(\hf + \delta
+ iv)|\int_{T^{1-\e}}^\infty x^{{1\over2}+\delta-\s-1}\d x\cr&
\ll t^{\e-1} + t^{2\e+{1\over2}-\s} \ll t^{2\e+{1\over2}-\s} \cr}
\leqno(4.13)
$$
if $1 < \s < {3\over2}$ and $\delta > 0$ is sufficiently small.
Finally we use (3.13), (4.9) (with $\s = \hf+\e$) and (4.13)
(with $\s = 1+\e$) to deduce that
$$
\Z_2(s) \ll_\e X^{{1\over2}-\s}t^\e + X^{1-\s}t^{\e-{1\over2}}
\ll t^{\e+{1\over2}-\s}
$$
for $X = t$. This finishes the proof of Corollary 4.

\smallskip
In concluding, note that (4.2) and (4.3) do not seem
to imply one another. It seems
even that the Riemann hypothesis (that all complex zeros of $\zeta(s)$
have real parts equal to $\hf$) does not imply (4.2). On the other hand,
(4.3) is a trivial consequence of the Lindel\"of hypothesis (4.2) (which
is a consequence of the  Riemann hypothesis; see [2] and [17]).

\bigskip
{\hh 5. The conjecture on exponential sums with Hecke series}
\bigskip
We pass now to a conjecture involving exponential sums with the Hecke series
$\H$, which will have applications to $\Z_2(s)$ and to $\zeta(s)$.
The Lindel\"of hypothesis (4.10) can be
recast (by using the approximate functional equation for $\zeta(s)$
and the Perron inversion formula [2, (A.10)])
in the form involving exponential sums, namely
$$
\sum_{N<n\le N'}n^{-it} \;\ll_\e\;N^{1\over2}t^\e
\qquad(N < N' \le 2N\ll t,\,
t > t_0 > 0).\leqno(5.1)
$$
It trivially implies the power moment estimates (see [2], [3] and [17]
for a comprehensive account)
$$
\int_0^{T}|\zt|^{2k}\d t \;\ll_{\e,k}  \; T^{1+\e}\qquad(k\in \NN),
\leqno(5.2)
$$
and in particular the eighth moment (4.3), namely the case $k=4$
of (5.2).
The connection with the Hecke series is that both estimates for $E_2(T)$ and
$\int_0^T|\zt|^8\d t$ can be made to depend on exponential sums with
the Hecke series $\H$ (see (5.8) and (5.12)-(5.14)).
\smallskip
One conjectures (this can be thought of as an analogue, in some sense,
of the Lindel\"of hypothesis) that
$$
H_j(\hf) \;\ll_\e\; \k_j^\e,\leqno(5.3)
$$
and more generally that
$$
H_j(\hf+it) \;\ll_\e\; (\k_j|t|)^\e.
$$
The conjecture (5.3) is, at the present state of knowledge involving
Hecke series, certainly out of reach. However, recently the author proved
in [7] that
$$
\sum_{K-G\le\k_j\le K+G} \a_j\H \;\ll_\e\; GK^{1+\e}\leqno(5.4)
$$
for
$$
K^{\e}  \;\le\; G \;\le \; K.\leqno(5.5)
$$
In view of the nonnegativity of $H_j(\hf)$ and $\a_j \gg_\e \k_j^{-\e}$,
this result implies that
$$
H_j(\hf) \;\ll_\e\; \k_j^{{1\over3}+\e}.\leqno(5.6)
$$
Of course, (5.6) is much weaker than the conjectural (5.3),
but nevertheless
it is the first published improvement over the trivial
$H_j(\hf) \ll \k_j^{{1\over2}}$. The results of (5.4)--(5.6) can be put in
the form
$$
\sum_{K-1\le\k_j\le K+1} \a_j\H \;\ll_\e\; K^{1+\e},\leqno(5.7)
$$
which is essentially best possible. However, when $\a_j\H \,(\ge 0)$
in (5.7) is
weighted by a suitable exponential factor, one expects additional
cancellation to take place, just like when instead of the sum
$$
\sum_{N<n\le N'}1 \= N' - N + O(1)\qquad(N< N'\le2N)
$$
we consider the sum in (5.1), which is weighted by the exponential factor
$\exp(-it\log n)$.

\medskip
In applications an exponential sum with $\a_j\H$ has occurred on
at least two important occassions. First, we have (proved by
Ivi\'c--Motohashi [9])
$$\eqalign{&
E_2(T) \ll \D\log^cT + \cr&
+ T^{1\over2}\sup_{\tau\asymp T}\left|
\sum_{\k_j\le T\D^{-1}\log T}\a_j\H\k_j^{-3/2}\exp\left(i\k_j
\log{\k_j\over\tau}-\left({\D\k_j\over T}\right)^2\right)\right|\cr}
\leqno(5.8)
$$
provided that $T^{1/2} \le \D \le T^{2/3}\log^cT$.

\medskip
The second application consists of the bound for the eighth moment
of $|\zt|$.
Note that in [10] it was proved that
$$
\int_T^{2T}|\zt|^8\d t \ll_\e T^{-2}\int_{T^{1/2}}^{T^{1+\e}}
t^2|{\cal G}(\hf+it)|^2\d t + T^{1+\e},\leqno(5.9)
$$
where for $\hf < a < 1$ one has
$$
{\cal G}(s) \;:=\; {1\over2\pi i}\int_{(a)}z_2(w)T^{w+s}\,{U(s,w)\over w}\d w,
\leqno(5.10)
$$
$$
U(s,w) \;:=\; \int_{1\over2}^{5\over2}\Phi(x)x^{s+w-1}\d x \ll
\min\left(1,\,{1\over|s+w|^A}\right)\leqno(5.11)
$$
for any given $A > 0$, where $\Phi(x) \in C^\infty$ is a nonnegative
 function supported
in $[\hf,\,{5\over2}]$ which equals unity in $[2,\,2]$.
From (5.9) and (5.11) we obtain, by shifting
appropriately the line of integration in the expression for
${\cal G}(s)$,
$$
\int_T^{2T}|\zt|^8\d t\ll_\e \int_{T^{1/2}}^{T^{1+\e}}t^2|I(T,t)|^2\d t
 + \Delta(T),
\leqno(5.12)
$$ where for the error term $\Delta(T)$ we hope to have
$\Delta(T) \ll_\e T^{1+\e}$, which is by no means easy to establish.
In (5.12) we have
$$
I(T,t) := \sum_{t-T^\e\le\k_j\le
t+T^\e}\alpha_jH_j^3(\hf){R_1(-\k_j)
\over{1\over2}-i\k_j}\int_{1/2}^{5/2}\Phi(x)(Tx)^{it-i\k_j}\d x,
\leqno(5.13)
$$
The function $R_1$ is closely connected to $R$ in (1.2), and we have
$$
R_1(y) := \sqrt{{\pi\over2}}{\Bigl(2^{-iy}{\G({1\over4} -{1\over2}
{i}y)
\over\G({1\over4} + {1\over2}{i}y)}\Bigr)}^3\G(2iy)\cosh(\pi y)
 \ll (1+|y|)^{-{1\over2}}\quad (y\in\RR).\leqno(5.14)
$$
If we  bound $I(T,t)$ by (5.7) we shall obtain
$$
I(T,t) \ll t^{-{3\over2}}\sum_{t-T^\e\le\k_j\le t+T^\e}\a_jH_j^3(\hf)
\ll t^{\e-{1\over2}}\qquad(T^{1/2} \le t \le T^{1+\e}),
$$
and consequently (5.12) gives the worse-than-trivial bound
$$
\int_0^T|\zt|^8\d t \ll_\e T^{2+\e}.
$$
The $x$--integral in (5.13) cannot help much because it is
practically non-oscillating. One does expect that massive cancellation
will be induced by $R_1(-\k_j)$. From Stirling's formula it follows that
$$
I(T,t) \;=\; O(t^{\e-{3\over2}}) \,+\leqno(5.15)
$$
$$
+ \,\pi(2t)^{-3/2}\sum_{t-T^\e\le\k_j\le
t+T^\e}\alpha_jH_j^3(\hf)\exp\left(i\k_j\log{\k_j\over4e}\right)
\int\limits_{1/2}^{5/2}\Phi(x)(Tx)^{it-i\k_j}\d x.
$$
The exponential sum involving
$\a_j\H$ becomes then essentially the same one as the sum in (5.8),
if the latter is split into short subsums.
Thus it seems reasonable to make the following

\bigskip
{\bf Conjecture 3}. {\it For $\tau^\delta \ll K \ll \tau^{1+\delta}
\,\;(0 < \delta < 1)\,$ we have}
$$
\sum_{K-1\le\k_j\le K+1}\a_j\H\exp\left(i\k_j\log{\k_j\over\tau}\right)
\;\ll_\e\;K^{{1\over2}+\e}.\leqno(5.16)
$$

\bigskip\noindent
Note that in (5.16) we are assuming a saving of $\sqrt{K}$ over
the known  bound (5.7) when there is no exponential
factor, which in a sense
corresponds to the saving required by the Lindel\"of hypothesis (4.10).
Very likely (5.16) is, if true, essentially best possible (see (7.6)).
Two consequences of (5.16) are the bounds (4.2) and (4.3),
similarly as if one assumes the conjecture (4.1).

\smallskip {\bf Corollary 5.} {\it If Conjecture 3 holds, then}
(4.2) {\it holds.}

\smallskip {\bf Corollary 6.} {\it If Conjecture 3 holds, then}
(4.3) {\it holds.}

\smallskip\no {\bf Proof of Corollary 5.}
To obtain (4.2)
we use (5.8), splitting the sum into $O(\log T)$ sums over $[K, K']
\;(K < K' \le 2K)$, and removing
the monotonic coefficients
$$
\k_j^{-3/2}\exp\left(-\left({\D\k_j\over T}\right)^2\right)
$$
by partial summation. Each of the sums over $[K, K']$ is further split
in $\ll K$ subsums over unit intervals, to which the conjecture (5.16) is
applied. The choice $\D = T^{1/2}$ gives then (4.2).

\smallskip\no {\bf Proof of Corollary 6.}
In a similar way
we use (5.16) to obtain that
$$
I(T,t)  \;\ll_\e\ t^{1+\e}\qquad(T^{1/2} \le t \le T^{1+\e}),
$$
which easily gives (4.3) by (5.12), provided one can prove rigorously
that $\Delta(T) \ll_\e T^{1+\e}$ in (5.12). This can be achieved by noting
that the main contribution to $\D(T)$ comes from the bound of the
portion of $\Z_{21}(s)$ (see (6.7)) coming from the discrete
spectrum at $\s = -\e$, when we shift the line of integration in
the relevant part of (5.10) to $a = -\e$. Then the relevant expression
will be an exponential sum with $\a_j\H$ to which (5.16) may be applied.
This will lead to $\Delta(T) \ll_\e T^{1+\e}$.

\smallskip
A possibility to treat the sum in (5.16)
is to use Motohashi's transformation formula (see [15, Lemma 3.8]
and the method of evaluating ${\cal C}(K,G)$ on p. 127) for
$$
\sum_{j=1}^\infty \a_jH_j^3(\hf)h(\k_j),\leqno(5.17)
$$
where $h(r)$ is an even function of exponential decay in a suitable
horizontal strip satisfying $h(\hf i) = 0$.
Instead of considering summation over the interval $[K-1,\,K+1]$, one could
consider summation over intervals of the form
$[K-G,\,K+G],\, K^\e \ll G \ll K^{1-\e},$ with the idea of
choosing $G$ suitably. A natural choice for the function $h$ is
$$\eqalign{
&h(r) = \cr&= (r^2 + {\txt{1\over4}})\cos\left(\hf r\log({r\over CT})^2
\right)\Bigl\{\exp\Bigl(-{({r-K\over G})}^2\Bigr)
+ \exp\Bigl(-{({r+K\over G})}^2\Bigr)\Bigr\}.\cr}
$$
This is a difficult problem, and even if
some progress with the sum (5.17) could be made,
this would not automatically imply any result concerning the sum (5.15),
where there is no Gaussian exponential factor.

\medskip
It may be also remarked that either of the Conjectures may be
used to derive a mean value result for Dirichlet polynomials. This is

\medskip {\bf Corollary 7}. {\it If  either Conjecture} 1,\enskip2 {\it or} 3
{\it holds, $1 \ll N \ll T^2 $ and
$a(n) \in\CC$ is arbitrary, then}
$$\eqalign{&
\int_0^T|\zt|^4\Bigl|\sum_{N<n\le N'\le2N}a(n)n^{it}\Bigr|^2\d t\cr&
\ll_\e T^\e\sum_{N<n\le N'\le2N}|a(n)|^2(T + T^{{1\over2}}N).\cr}\leqno(5.18)
$$
Corollary 7 may be compared to the result of Deshouillers--Iwaniec [1],
 who had $(T + T^{1/2}N^2 + T^{3/4}N^{5/4})$
 as the factor on the right-hand side of (5.18).
Assuming the Selberg conjecture
that the smallest positive eigenvalue of the non-Euclidean Laplacian
for Hecke congruence subgroups is $\ge 1/4$, then the
Deshouillers--Iwaniec proof shows that the term
$T^{3/4}N^{5/4}$ in the above factor may be discarded.
This  is the limit of the Deshouillers--Iwaniec method.
N. Watt in [18] showed that  the left-hand side of (5.18) is
$\ll_{\e} T^{\e}(T + T^{1/2}N^2)N\max_{N<n\le 2N}|a_n|^2$, which
improves the Deshouillers--Iwaniec bound if $a_n \ll_\e n^\e$.

\bigskip
{\hh 6. A new bound for $\Z_2(s)$}
\bigskip

We have seen that both Conjecture 2 and Conjecture 3 imply the important
bounds (4.2) and (4.3) in zeta-function theory. Thus the question
naturally arises: is there any connection between Conjecture 2
and Conjecture 3?
Does one of them imply the other? They appear to be both of the
same level of difficulty, and it will be shown now heuristically
that Conjecture 3 implies Conjecture 2. We start from (3.13) and replace
the left side of ${\cal D}$  by the line $\Re{\rm e}\,w = -\e$.
Then we expect that the major contribution should come from the poles
at $w = \hf \pm i\k_j-s$. In view of the gamma-factor it transpires then
that we obtain the relevant sum
$$
\sum_{|\k_j-t|\le t^\e}\a_j\H X^{{1\over2} +
i\k_j-s}R_1(\k_j)\G(\hf +i\k_j-s),
\leqno(6.1)
$$
where $R_1$ is given by (5.14).
If we disregard the exponential factor which will come
from $R_1(\k_j)$  and use only use (5.7),
then from (3.13) and (3.14) (with $X =t$) we obtain the bound
$$
\Z(\s + it) \;\ll_\e\; t^{1-\s+\e}\qquad(\hf < \s < 1),\leqno(6.2)
$$
which improves the bound obtained in [10], where the exponent of $t$
was $2-2\s+\e$. However, if we use the Conjecture 3,
 then the sum in (6.1) will be $\ll_\e t^\e X^{{1\over2}-\s}$,
and Conjecture 2 follows. We shall derive now rigorously (6.2) and prove
the following

\medskip
{\bf THEOREM 3}. {\it If $s=\s+it$ is well separated from the poles of
$\Z_2(s)$, then for $0 < \s < 1, \,t \ge t_0 > 0$ we have}
$$
\Z_2(\s + it) \;\ll_\e\; t^{1-\s+\e}.\leqno(6.3)
$$

\bigskip\no
{\bf Proof.} From the bounds [10, (4.13) and (4.29)] we have,
under the above hypotheses,
$$
\Z_2(\s + it) \;\ll_\e\; t^{{1-\s\over1-\xi}+\e} +
\sum_{|t-\k_j|\le t^\e}\a_j\H t^{{\xi-2\s\over2-2\xi}+\e},\leqno(6.4)
$$
provided that
$$
{\txt{1\over3}} \le\xi\le\hf.\leqno(6.5)
$$
Take first $\s = \hf + \delta,\,\xi = {1\over3}$. Then from (6.4)
and (5.7) we obtain
$$
\Z_2(\s + it) \;\ll_\e\; t^{{3\over2}(1-\s)+\e},
$$
which is weaker than (6.3). The desired bound (6.3) would clearly
follow if $\xi = \e$ is a permissible value in (6.4).
To ascertain this fact there are two ways to proceed. Y. Motohashi
[13], [15] established the spectral decomposition of the function
$$
\psi(T) \;= \; \psi(T,\xi) \;:= \;
{1\over\sqrt{\pi}T^{\xi}}\int_{-\infty}^{\infty}\z {(T+t)}
\exp(-(t/T^{\xi})^2)\d t,\leqno(6.6)
$$
where $0 < \xi < 1$ is a constant.
In [10] this function was used
as a substitute for  $\z T$ in dealing with $\Z_2(s)$ for
$\s \le \hf$.  Namely we put
$$
\eqalign{ {\cal Z}_2(s)&=\int_1^{\infty }\psi(T)T^{-s}\d T + \int_1^{\infty
} (
\z T -\psi(T))T^{-s}\d T\cr &={\cal Z}_{21}(s)+{\cal Z}_{22}(s),
\cr}\leqno(6.7)
$$
say. It turns out that the integral ${\cal Z}_{22}(s)$ will converge well,
and the main difficulties are inherent in ${\cal Z}_{21}(s)$. The key
r\^ole in the spectral decomposition of $\psi(T)$ is played by the
function  ($r$ is real)
$$
\eqalign{
\Xi (ir;T,T^{\xi })\;=\;&{\G^2 ({1\over2}+ir)
\over{\G (1+2ir)}}\int_0^{\infty }
(1+y)^{-1/2+iT}y^{-1/2+ir}\times\cr
& \exp\left(-\txt{1\over4} T^{2\xi}\log ^2(1+y)\right )
F(\hf +ir, \hf + ir;1+2ir;-y) \d y,\cr}
$$
where $F$ is the hypergeometric function. For this function Motohashi
[15, (5.1.39)--(5.1.41)] obtains an asymptotic formula, where he essentially
has the condition (6.5). The first approach is to go carefully through
Motohashi's proof of the asymptotic formula for $\Xi (ir;T,T^{\xi })$,
and try to relax the condition (6.5) to $0 < \xi \le \hf$. This can be
done, but the analysis is rather long and technical. We shall adopt here
another approach. This consists of going through the proof of (6.4) in [10],
and making appropriate modifications.
The condition (6.5) is actually used there in the estimation of the
contribution coming from the saddle point $z_0$, the root of
$F(z) = 0$, where
$$
F(z) = F(z;r,T) = -r\log z + T\log(1 +{z\over T}) + 2r\log\left
(1 + \sqrt{1 + {z\over T}}\,\right),
$$
so that
$$
F'(z) = -{r\over z} + {T\over T+z}
+{r\over T\left(\sqrt{1 + {z\over T}}
+ 1 + {z\over T}\right)}.
$$
We obtain
$$
{z_0\over r} = \left(1 + {z_0\over T}\right)
{\left(1 + {r/T\over1+(1 + {z_0\over T})^{-1/2}}\right)}^{-1},
$$
and since it was shown in [10] that we have the conditions
$$
T \;\ge\; T(r) \;:=\; r^{1\over1-\xi}\log^Dr,
\qquad|r - t| \le t^\e,\leqno(6.8)
$$
we find by iteration that
$$
{z_0\over r} \= 1  + {r\over2T} + {r^2\over8T^2} + O\left({r^3\over T^3}
\right),\leqno(6.9)
$$
where the $O$--term in (6.9) admits
asymptotic expansion in term of powers of $r/T$.
The crucial term in the contribution of $z_0$ is $G(z_0)$, where
$$\eqalign{
G(z) &\;:=\;
i\log(1 + {z\over T})
+ {z\over T^2}{(\hf-iT)\over 1 + {z\over T}}\cr&
+ {\left(1 + \sqrt{1 + {z\over T}}\,\right)}^{-1}{(1-2ir)z
\over 2T^2\sqrt{1 + {z\over T}}}
-\hf\xi T^{2\xi-1}\log^2(1 + {z\over T})\cr&
+ \hf T^{2\xi}\log(1 + {z\over T}){z\over T^2}.\cr}
$$
so that
$$
G(z_0) = {iz_0\over2T^2}(z_0-r)-{2iz_0^3\over3T^3} + O\left({r\over T^2}
\right) + O(r^2T^{2\xi-3}). \leqno(6.10)
$$
In [10] (cf. page 332, line 7) we estimated trivially $G(z_0)$ as
$$
G(z_0) \ll r^3T^{-3} + rT^{-2} + r^2T^{2\xi-3} \ll  rT^{-2} +
r^{2+\e}T^{2\xi-3},\leqno(6.11)
$$
since in view of (6.5) and (6.8) we have
$$
r^3T^{-3} \le r^2\log^Cr\cdot T^{-2-\xi} \ll_\e r^{2+\e}T^{2\xi-3}.
$$
Then by trivial estimation the total contribution coming from the function
$\Z_{21}(s)$ is bounded by the second term on the right-hand side of
(6.4). However we can deal with the first two terms on the right-hand
side of (6.10) as follows. First note that by (6.8) and (6.9)
$$
{iz_0\over2T^2}(z_0-r)-{2iz_0^3\over3T^3} \;\asymp\;{r^3\over T^3}.
$$
The crucial portion of the integral
$$
X_r^*(s) \;:=\; \int_{T(r)}^{\infty}\Xi (-ir;T,T^{\xi})T^{-s}\d T
$$
which is to be bounded is of the form
$$
(\hf+ir-s)^{-1}r^{-1/2}\int_{T(r)}^\infty K(r,T)e^{iF(z_0)}T^
{{1\over2}+ir-s}\d T,\leqno(6.12)
$$
where $K(r,T) \asymp (r/T)^3$ and (6.8) is assumed to hold. A
calculation shows that
$$
{\partial F(z_0)\over\partial T} \;\gg\; {r^2\over T^2}\leqno(6.13)
$$
holds in the relevant range for $r$ and $T$. Write
$$
\int_{T(r)}^{\infty} = \int_{T(r)}^{\infty}K(r,T)T^{{1\over2}-\s}\cdot
e^{iF(z_0)+i(r-t)\log T}\d T,
$$
set
$$
H(r,T) \= F(z_0) + (r-t)\log T
$$
and observe that, since (6.8) and (6.13) hold,
$$
{\partial H(r,T)\over\partial T} \;\gg\; {r^2\over T^2}\qquad({\rm
for}\; t \ge T^{{1\over2}+\e}).
$$
Then we have by the first derivative test (cf. [2, Lemma 2.1]),
for $t \ge T^{{1\over2}+\e}$, that the total contribution of
$r = \k_j$ does not exceed the second term on the right-hand
side of (6.4). In case when $t < T^{{1\over2}+\e}$, that is,
$T > t^{2-\e}$, we obtain by trivial estimation that the contribution
of the integral in question is
$$\eqalign{&
\ll_\e |\hf+ir-s|^{-1}r^{-1/2}\int_{t^{2-\e}}^\infty {r^3\over T^3}
T^{{1\over2}-\s}\d T\cr&
\ll_\e |\hf+ir-s|^{-1}r^{5/2}t^{(2-\e)(-{3\over2}-\s)} \ll_\e
|\hf+ir-s|^{-1}t^{-{1\over2}-2\s+\e}.\cr}
$$
If $s = \s +it$ is well separated from the poles of $\Z_2(s)$,
then by (6.8) and (5.7) the total contribution of the portion in
question is
$$
\ll_\e \sum_{|t-\k_j|\le t^\e}\a_j\H\cdot t^{-{1\over2}-2\s+\e} \ll_\e
t^{{1\over2}-2\s+\e} \ll_\e t^{1-\s+\e}\qquad(\s \ge -\hf),
$$
which is negligible. Since there is now no restriction on $\xi$ except
the initial one that $0 < \xi \le \hf$, we obtain (6.3) from (6.4)
with $\xi = \e$, as asserted. We note that one can extend the validity
of (6.3) to the half-plane $\s > -\hf$, since in [10] it was shown that
(6.4) holds for
$$
\s \;>\; \max(\,2\xi - 1,\,-\hf\,).
$$
In the other direction (6.3) holds at least for $1 \le \s \le 2$,
which follows from (3.24).

\smallskip
We note that the bound
$$
\Z_2(\hf + \e + it) \ll \sqrt{t},\leqno(6.14)
$$
which follows
from (6.3), implies the mean square bound (3.15),  therefore also
(3.16) and (3.19). More generally, if we assume that
$$
{\cal Z}_2(\hf + \e +it) \;\ll\; t^\rho
$$
holds with some $0 < \rho \le \hf$, then from (4.6) we obtain
$$
E_2(T) \;\ll_\e\; T^{{2\rho+1\over2+2\rho}+\e}.\leqno(6.15)
$$
The bound (6.15) clearly shows
that any improvement of (6.3) would have far-reaching
consequences in zeta-function theory. It also transpires that
essentially progress on bounds for $E_2(T)$ and the eighth moment
of $|\zt|$ follows from new bounds for the exponential sums appearing
in (5.16).

\bigskip\bigskip
{\hh 7. The function $\Z_k(s)$}
\bigskip

The function ${\cal Z}_2(s)$ is a special case of the function
$$
{\cal Z}_k(s) := \int_1^\infty|\zx|^{2k}x^{-s}\d x\qquad(s
 = \s + it;\; \s,t \in\RR, \, k\in\NN),
\leqno(7.1)
$$
introduced in [10], where the cases
$k = 1,2$ were extensively investigated.
One of the possible applications of ${\cal Z}_k(s)$ consists of the
following. If $F(s)$ is the Mellin transform of $f(x)$,
then by (2.1) one formally obtains, for suitable $c > 1$,
$$\eqalign{
\int\limits_1^\infty f\left({x\over T}\right)|\zx|^{2k}\d x &=
\int\limits_1^\infty{1\over2\pi i}
\int\limits_{(c)}F(s)\left({T\over x}\right)^s\d s|\zx|^{2k}\d x\cr&
= {1\over2\pi i}\int_{(c)}F(s)T^s{\cal Z}_k(s)\d s.\cr}\leqno(7.2)
$$
If $f(x) \in C^\infty$ is a nonnegative function of compact support
such that $f(x) = 1$ for $1 \le x \le 2$, then $F(s)$ is entire of
fast decay, and (7.2) (with $c = 1+\e$) yields a weak form of the
$2k$--th moment for $|\zt|$, namely
$$
\int_0^T\vert\zeta({\txt{1\over 2}} + it)\vert^{2k}\d  t\;
\ll_{k,\e} \; T^{1+\e},\leqno(7.3)
$$
provided that ${\cal Z}_k(s)$ has analytic continuation to the half-plane
$\s > 1$, where it is regular and of polynomial growth in $|t|$.
Conversely, if (7.3) holds, then integrating by parts the right-hand
side of (7.1) it is seen that ${\cal Z}_k(s)$  is regular for $\s > 1$
and in this half-plane satisfies ${\cal Z}_k(s) \ll_\s |t|$. Thus the
$2k$--th moment estimate (7.3) has an equivalent formulation
in terms of the analytic behaviour of ${\cal Z}_k(s)$. Moreover (7.3)
is equivalent (for example, by Lemma 7.1 of [3]) to the
Lindel\"of Hypothesis. Hence it can be said that
the Lindel\"of hypothesis is equivalent to the following statement:
given $\e > 0$, for every $k \in \NN$ the function ${\cal Z}_k(s)$
can be analytically continued to the half-plane $\s \ge 1 + \e$,
where it is of polynomial growth in $|t|$.

\medskip
In [10] it was shown that
$$
\int_1^T|\Z_1(\s + it)|^2\d t \;\ll_\e\; T^{2-2\s+\e}\qquad(\hf \le \s \le
1).\leqno(7.4)
$$
By arguments analogous to the ones used in the proof of Theorem 1 it can
be shown that (7.4) is essentially best possible, namely that the
integral in question is $\gg_\e T^{2-2\s-\e}$. Another, quicker proof
of (7.4) follows from the method of proof of Theorem 3. This in fact will
yield even
$$
\int_T^{2T}|\Z_1(\s + it)|^2\d t \;\ll_\e\; T^{2-2\s+\e} +T^{-1}\qquad
(\s \ge \hf).\leqno(7.5)
$$
The bound (7.4) is the analogue of Conjecture 1 (cf. (4.1)) for $\Z_1(s)$.
However, the analogue of Conjecture 2 (cf. (4.9)), namely
$$
\Z_1(\s + it) \;\ll_\e\;|t|^\e\qquad(\s > \hf),
$$
is a difficult problem which is certainly out of reach at present.

\medskip
One can directly estimate the $2m$-th moment of $|\zx|$ by ${\cal Z}_m(s)$
as follows. Let $m = k + \ell,\, k,\ell \in \NN$. Then by Lemma 1
we obtain, provided $c$ and $d$ are sufficiently large
$$\eqalign{
&\int_{T/2}^{5T/2}\f(x)|\zx|^{2m}\d x \cr&
= {1\over2\pi i}\int_{(c)}\Z_k(s)\int_{T/2}^{5T/2}\f(x)
|\zx|^{2\ell}x^{s-1}\d x\,\d s
\cr&
= \left({1\over2\pi i}\right)^2\int_{(c)}\Z_k(s)\int_{(d)}\Z_\ell(w)
\int_{T/2}^{5T/2}\f(x)x^{s+w-2}\d x\,\d s\,\d w.
\cr}\leqno(7.6)
$$
Integrating by parts sufficiently many times the integral over $x$
it is seen that the non-trivial contribution in (7.6) comes from
$w,\,s$ satisfying $|v+t| \le |t|^\e$, where $s = \s+it,\,w = u + iv$.
Then from (7.6) we obtain
$$
 \int_T^{2T}|\zx|^{2m}\d x \ll T^{c+d-1},
$$
provided that $c, d$ are chosen in such a way that the triple
integral in (7.6) converges. This procedure shows that mean values
of ${\cal Z}_k(s)$ are relevant not only in the range $\s < 1$, but
in the range $\s \ge 1$ as well (see (4.6)). A way to deal with
${\cal Z}_m(s)$ is to observe that, by Lemma 2, one has for
sufficiently large $c$,
$$
{\cal Z}_m(s) \= {1\over2\pi i}\int_{(c)}
{\cal Z}_\ell(w){\cal Z}_k(s+1-w)\d w.\leqno(7.7)
$$
Thus (7.7) is a sort of a recurrent relation that permits one to
deduce information on ${\cal Z}_m(s)$ from ${\cal Z}_k(s)$
and ${\cal Z}_\ell(s)$ if $m = k + \ell$.

\smallskip
A possible application of (7.6) is to show that if one assumes
Conjecture 3, then it follows that it is essentially optimal.
Namely we expect that there exists a constant $C > 0$ such that,
for given $\delta,\,\e >0$, we have
$$
\sup_{T^{1-\delta}\le K\le T^{1+\delta}}
\left|\sum_{K-1\le\k_j\le K+1}\a_j\H\exp\left(i\k_j\log{\k_j\over CT}
\right)\right| \;\gg_{\delta,\e}\; T^{{1\over2}-\e}.\leqno(7.8)
$$
This can be shown heuristically as follows. Take $m =4, k = \ell = 2$
in (7.6), so that the left-hand side is $\gg T\log^{16}T$. From the
first identity in (7.6), by using repeated integration by parts
and the fact (see Section 3) that $Z^4(t)$ and its derivatives are
bounded in mean by log-powers, it follows that the $s$--integral
in (7.6) can be truncated at $|\I s| = T^{1+\delta}$ with error
which is $\ll 1$ when we take $c = \hf + \e$. In the $w$--integral
the relevant portion (coming from the discrete spectrum in $\Z_{21}(s)$
in (6.7)) is obtained by taking $a = -\e$. By using repeated integration
by parts in the $x$--integral it transpires that only the terms $|t + \k_j|
\le t^\e$ are relevant. Therefore we obtain
$$\eqalign{&
T\log^{16}T \ll \int\limits_{-T^{1+\delta}}^{T^{1+\delta}}
|\Z_2(\hf+\e+it)|\times\cr&\times\left|\sum_{|\k_j+t|\le t^\e}
\a_j\H R_1(-\k_j)\int\limits_{T/2}^{5T/2}\f(x)x^{\e-1+t+i\k_j}\d x
\right|\d t.\cr}\leqno(7.9)
$$
Conjecture 3, as was seen in Section 6, implies heuristically the truth
of Conjecture 2. Thus we may truncate the $t$-integral in (7.9)
at $|t| = T^{1-\delta}$ with an error which will be $\ll T$ if
$\delta,\e$ are sufficiently small. Then we use (5.14) and Stirling's
formula to obtain an asymptotic formula for the function $R_1$.
Finally we divide the range of summation over the spectrum into
$O(T^\e)$ subsums with the range of summation $K-1\le\k_j\le K+1$,
interchange summation and integration over $x$ and take the
suprema over $x$ and $K$ to obtain (7.8) after trivial estimation.

\bigskip\bigskip
\centerline{\hh REFERENCES}
\bigskip

\item {[1]} J.-M. Deshouillers and H. Iwaniec,  Power mean values of
the Riemann zeta-function, {\it Mathematika} {\bf29}(1982), 202-212.

\item {[2]} A. Ivi\'c,  The Riemann zeta-function, {\it John Wiley and
Sons}, New York, 1985.

\item {[3]} A. Ivi\'c,  Mean values of the Riemann zeta-function, LN's
{\bf 82}, {\it Tata Institute of Fundamental Research}, Bombay,  1991
(distr. by Springer Verlag, Berlin etc.).

\item{[4]}  A. Ivi\'c,  On the fourth moment of the Riemann
zeta-function, {\it Publs. Inst. Math. (Belgrade)} {\bf 57(71)}
(1995), 101-110.

\item{[5]} A. Ivi\'c,  The Mellin transform and the Riemann
zeta-function,  {\it Proceedings of the Conference on Elementary and
Analytic Number Theory  (Vienna, July 18-20, 1996)},  Universit\"at
Wien \& Universit\"at f\"ur Bodenkultur, Eds. W.G. Nowak and J.
Schoi{\ss}engeier, Vienna 1996, 112-127.

\item{[6]} A. Ivi\'c, On the error term for the fourth moment of the
Riemann zeta-function, {\it J. London Math. Soc.}
{\bf60}(2)(1999), 21-32.

\item{[7]} A. Ivi\'c, On sums of Hecke series in short intervals,
to appear in J. de Th\'eorie des Nombres Bordeaux.

\item{ [8]} A. Ivi\'c and Y. Motohashi,  The mean square of the error
term for the fourth moment of the zeta-function,  {\it Proc. London
Math. Soc.} (3){\bf 66}(1994), 309-329.

\item {[9]} A. Ivi\'c and Y. Motohashi,  The fourth moment of the
Riemann zeta-function, {\it J. Number Theory} {\bf 51}(1995), 16-45.

\item {[10]} A. Ivi\'c, M. Jutila and Y. Motohashi, The Mellin
transform of powers of  the Riemann zeta-function,  Acta Arith.
{\bf95}(2000), 305-342.

\item {[11]} A.A. Karacuba and S.M. Voronin, The Riemann zeta-function,
Walter de Gruyter, Berlin--New York, 1992.

\item {[12]} S. Katok and P. Sarnak, Heegner points, cycles and Maass
forms, {\it Israel J. Math.} {\bf84}(1993), 193-227.

\item{ [13]} Y. Motohashi,   An explicit formula for the fourth power
mean of the Riemann zeta-function, {\it Acta Math. }{\bf 170}(1993),
181-220.

\item {[14]} Y. Motohashi,  A relation  between the Riemann
zeta-function and the hyperbolic Laplacian, {\it Annali Scuola Norm.
Sup. Pisa, Cl. Sci. IV ser.} {\bf 22}(1995), 299-313.

\item {[15]} Y. Motohashi,  Spectral theory of the Riemann
zeta-function, {\it Cambridge University Press}, Cambridge, 1997.

\item {[16]} E.C. Titchmarsh, Introduction to the Theory of Fourier
Integrals, {\it Clarendon Press}, Oxford, 1948.

\item {[17]} E.C. Titchmarsh,  The Theory of the Riemann Zeta-Function
(2nd ed.),
 {\it Clarendon  Press}, Oxford, 1986.

\item{[18]} N. Watt, Kloosterman sums and a mean value
for Dirichlet polynomials, {\it J. Number Theory} {\bf53}(1995),
179-210.

\bigskip
\parindent=0pt

\cc
Aleksandar Ivi\'c \par
Katedra Matematike RGF-a
Universiteta u Beogradu\par
Dju\v sina 7, 11000 Beograd,
Serbia (Yugoslavia)\par
{\sevenbf e-mail: aivic@matf.bg.ac.yu,
aivic@rgf.bg.ac.yu}

\bye